\documentclass[preprint, 12pt]{myelsarticle}
\usepackage{amssymb}
\usepackage[displaymath]{lineno}
\usepackage{graphicx}
\usepackage{caption}
\usepackage{subcaption}
\usepackage{sectsty}
\usepackage{natbib}
\usepackage{amsmath}
\journal{}
\newtheorem{theorem}{Theorem}[section]

\newtheorem{proposition}[theorem]{Proposition}
\newtheorem{corollary}[theorem]{Corollary}

\newtheorem{example}[theorem]{Example}

\newtheorem{remark}[theorem]{Remark}
\newproof{proof}{Proof}

\newcommand{\C}{{\mathbb C}}
\newcommand{\R}{{\mathbb R}}
\newcommand{\N}{{\mathbb N}}

\newcommand{\Z}{{\mathbb Z}}

\newcommand{\zp}{\lbrace}
\newcommand{\zz}{\rbrace}

\newcommand{\D}{{\rm\Delta}}
\newcommand{\bD}{\partial\D}
\newcommand{\bdr}{\partial}

\begin{document}
\begin{frontmatter}
\title{Mixed Riemann-Hilbert boundary value problem\\ with simply connected fibers}
\author{Miran \v{C}erne\fnref{thanks}}
\fntext[thanks]{
The author was supported in part by grants
{\it Analiza in geometrija} P1-0291, 
{\it Kompleksna in geometrijska analiza} J1-3005, 
{\it Holomorfne parcialne diferencialne relacije} N1-0237 and
{\it Nelinearni valovi in spektralna teorija} N1-0137
from ARRS, Republic of Slovenia.
}
\ead{miran.cerne@fmf.uni-lj.si}
\affiliation{organization={Faculty of Mathematics and Physics, University of Ljubljana, and Institute of Mathematics, Physics and Mechanics},
            addressline={Jadranska 19}, 
            city={Ljubljana},
            postcode={1\,111},
            country={Slovenia}}
\begin{abstract}
We study the existence of solutions of 
mixed Riemann-Hilbert or Cherepanov boundary value problem with simply connected fibers on the unit disk $\D$.
Let $L$ be a closed arc on $\bdr\D$ with the end points $\omega_{-1}, \omega_1$ and let $a$ be a smooth function on $L$ with no zeros. Let 
$\zp\gamma_{\xi}\zz_{\xi\in\bD\setminus\mathring{L}}$ be a smooth family of smooth Jordan curves in $\C$ which all contain point $0$
in their interiors and such that $\gamma_{\omega_{-1}}$, $\gamma_{\omega_{1}}$ are strongly starshaped with respect to $0$.
Then under condition that for each 
$w\in \gamma_{\omega_{\pm 1}}$ the angle between $w$ and the normal to $\gamma_{\omega_{\pm 1}}$
at $w$ is less than $\frac{\pi}{10}$,
there exists a H\"{o}lder continuous function $f$ on $\overline{\D}$, holomorphic on $\D$, such that 
\begin{linenomath*}
$${\rm Re}(\overline{a(\xi)} f(\xi)) = 0 {\rm\ on\  } L\ \ \ \ \ {\rm and}\ \ \ \ \ 
f(\xi)\in\gamma_{\xi} {\rm\ on\   } \bdr\D\setminus\mathring{L}.$$
\end{linenomath*}
\end{abstract}

\end{frontmatter}

\section{Introduction}
Let $\D = \zp z\in\C; |z|<1\zz$ be the open unit disc in the complex plane $\C$ and let
$\bD = \zp \xi\in\C; |\xi|=1\zz$ be the unit circle. 
Let $L$ be a closed arc on $\bdr\D$, let $\mathring{L}$ denote its interior with respect to $\bdr\D$, and let $a : L\rightarrow\bD$
be a smooth function.

Recall that the interior ${\rm Int}(\gamma)$ of a Jordan curve $\gamma\subseteq\C$
is the bounded component of $\C\setminus\gamma$. We orient $\gamma$ positively with respect to  ${\rm Int}(\gamma)$.
Jordan curve $\gamma\subset\C$ is {\it starshaped} with respect to $0$, if for any point $w$
in the interior of $\gamma$ the line segment which connects points $0$ and $w$ lies in
the  interior of $\gamma$, 
and it is {\it strongly starshaped} with respect to $0$, \cite{Han-Her-Mar-Mos}, if  there exists 
a positive continuous function $R$ on the unit circle such that
\begin{linenomath*}
\begin{equation}\label{stronglystar}
\gamma=\left\{ w\in\C; |w|= R\left(\frac{w}{|w|}\right)\right\}
\end{equation}
\end{linenomath*}
and 
\begin{linenomath*}
\begin{equation}\label{stronglystarint}
{\rm Int}(\gamma)=\left\{ w\in\C\setminus\{0\}; |w|< R\left(\frac{w}{|w|}\right)\right\}\cup\{0\}.
\end{equation}
\end{linenomath*}

Let $\zp\gamma_{\xi}\zz_{\xi\in\bD\setminus\mathring{L}}$ be a smooth family of smooth Jordan curves in $\C$ which all contain point $0$
in their interiors.
In this paper we study the existence and properties of {\it holomorphic} solutions of the nonlinear mixed Riemann-Hilbert problem, that is, the Che\-re\-pa\-nov boundary value problem with simply connected fibers. The problem asks for a continuous function
$f$ on $\overline{\D}$, holomorphic on $\D$, such that
\begin{linenomath*}
\begin{equation}\label{eq1L}
{\rm Re}(\overline{a(\xi)} f(\xi)) = 0 {\rm\ for\ } \xi\in L  
\end{equation}
\end{linenomath*}
and
\begin{linenomath*}
\begin{equation}\label{eq1K}
f(\xi)\in\gamma_{\xi} {\rm\ for\   } \xi\in\bdr\D\setminus\mathring{L}.
\end{equation}
\end{linenomath*}
That is, $f$ solves a linear Riemann-Hilbert problem on $L$ and a nonlinear Riemann-Hilbert
problem with simply connected fibers on $\bdr\D\setminus\mathring{L}$. See also \cite{Che2, Che1, Mit-Rog, Obn-Zul}.

The problem with circular fibers $\gamma_{\xi}$ and $L$ a finite union of disjoint arcs was considered by Obnosov and Zulkarnyaev in \cite{Obn-Zul}, and by the author in \cite{Cer2}.
The structure of the family of solutions of problem (\ref{eq1L}-\ref{eq1K}) is well known in the cases where
either $L=\bdr\D$ or $L=\emptyset$. If $L=\bdr\D$, we consider
a homogeneous linear Riemann-Hilbert problem. In this case the essential information on the problem is given by the winding number $W(a)$ of function $a$. It is well known \cite{Kop, Weg1, Weg2} that if the winding number $W(a)$ is nonnegative,
the space of solutions of (\ref{eq1L}) is a vector subspace of $A^{\alpha}(\D)$, $0<\alpha<1$, of real dimension $2W(a)+1$. 

\begin{remark}
{\rm The linear Riemann-Hilbert problem can also be considered in the case of a nonorientable line bundle over $\bdr\D$, that is, in the case where at some point $\xi_0\in\bdr\D$
we have $a(\xi_0^-) = - a(\xi_0^+)$. 
Then the winding number of function
$a^2$ or the Maslov index of the problem is an odd integer. In this case it holds that if $W(a^2)\ge -1$, or, with a little bit of abuse of notation, if $W(a)\ge -\frac{1}{2}$, then 
the space of solutions of (\ref{eq1L}) is a vector subspace of $A^{\alpha}(\D)$ of real dimension $2W(a)+1$, see \cite{Cer0, Cer1, Oh, Weg2}.}
\end{remark}

If $L$ is empty, we have a nonlinear
Riemann-Hilbert problem with smooth simply connected fibers which all contain $0$ in their interiors. This problem was considered and solved in \cite{For, Sni, Weg1, Weg2}. In particular, it was proved that
the family of solutions with exactly $m$ zeros on $\D$,  $m\in\N\cup\{0\}$, forms a manifold in
space $A^{\alpha}(\D)$ of dimension $2m+1$, and this manifold is compact if and only if $m=0$. We assume from now on that neither $L=\emptyset$ nor $L=\bdr\D$.

\begin{theorem}\label{th1}
Let $k\ge 3$. Let  $a : L\rightarrow\C\setminus\{0\}$ be a $C^{k+1}$ function and
let $\zp\gamma_{\xi}\zz_{\xi\in\bD\setminus\mathring{L}}$ be a $C^{k}$ family of Jordan curves in $\C$ which all contain 
point $0$ in their interiors. Let $\omega_{1}$ and $\omega_{-1}$ be the first and the last point of arc $L$ with respect to the positive orientation of $\bdr\D$. Let Jordan curves $\gamma_{\omega_j}$, $j=\pm 1$, be strongly starshaped with respect to $0$ and such that for each 
$w\in \gamma_{\omega_{j}}$ the angle between $w$ and the outer normal to $\gamma_{\omega_{j}}$
at $w$ is less than $\frac{\pi}{10}$.
Let $w_j$, $j=\pm 1$, be the intersection of $\gamma_{\omega_j}$ and the line ${\rm Re}(\overline{a(\omega_j)} w)=0$ of the form $\lambda (-i a(\omega_j))$, $\lambda>0$,
and let $\pi\beta_j$ be the oriented angle of intersection of the line ${\rm Re}(\overline{a(\omega_j)}w)=0$ 
with the fiber $\gamma_{\omega_j}$ at point $w_j$,  where $\beta_1\in(0,1)$ and $\beta_{-1}\in(-1,0)$.  
Let
\begin{linenomath*}
\begin{equation}
0< \beta < \min\{\beta_1, 1-\beta_1, |\beta_{-1}|, 1-|\beta_{-1}|\}.
\end{equation}
\end{linenomath*}
Then there exists a unique $f\in A^{\beta}(\D)$  with no zeros on $\D$ which solves
(\ref{eq1L}-\ref{eq1K}) for which $f(\omega_1)=w_1$ and $f(\omega_{-1})=w_{-1}$.
\end{theorem}

\begin{remark}
{\rm Here $\beta_1> 0$, if the tangent vector $ -i a(\omega_1)$ to ${\rm Re}(\overline{a(\omega_1)}w)=0$ 
is rotated counterclockwise by angle $\pi\beta_1$ to get a positive tangent vector to $\gamma_{\omega_1}$ at point $w_1$, and $\beta_{-1}<0$, if
a positive tangent vector to $\gamma_{\omega_{-1}}$ at $w_{-1}$
is rotated clockwise by angle $\pi|\beta_1|$ to get tangent vector $ -i a(\omega_{-1})$
to ${\rm Re}(\overline{a(\omega_{-1})}w)=0$.
}
\end{remark}

\begin{remark}
{\rm Observe that conditions in Theorem \ref{th1} imply $||\beta_j| - \frac{1}{2}|<\frac{1}{10}$, $j=\pm 1$, and hence one could choose $\beta=\frac{2}{5}$.}
\end{remark}

\begin{remark}
{\rm In the cases considered in \cite{Cer2, Obn-Zul} all boundary curves were circles with center at point $0$. Hence $|\beta_j| = \frac{1}{2}$, $j=\pm 1$, and the maximal regularity we got was $\beta<\frac{1}{2}$.}
\end{remark}

\begin{corollary}\label{cor1}
Let $a_1,\dots, a_n\in\D$ be a finite set of points with given multiplicities.
Then under the assumptions of Theorem \ref{th1}
there exists $\beta\in (0,1)$ and $f\in A^{\beta}(\D)$ which has zeros exactly at points $a_1,\dots, a_n\in\D$ with the given multiplicites
and which solves (\ref{eq1L}-\ref{eq1K}).
\end{corollary}

\section{Function spaces, Hilbert transform and\\ defining functions}

Let $0<\alpha<1$ and let $G\subset\C$ be a compact subset. We denote by $C^{\alpha}(G)$ the algebra over $\C$
of H\"{o}lder continuous complex functions on $G$ and by $C_{\R}^{\alpha}(G)$ the 
algebra over $\R$ of real H\"{o}lder continuous functions on $G$. Using the norm
\begin{linenomath*}
\begin{equation}
\| f\|_{\alpha} = \max_{z\in G} |f(z)|  + \sup_{z,w\in G, z\ne w}\frac{|f(z)-f(w)|}{|z-w|^{\alpha}}
\end{equation}
\end{linenomath*}
the algebras $C^{\alpha}(G)$ and $C_{\R}^{\alpha}(G)$ become Banach algebras. 
For $G=\overline{\D}$ or $G=\bdr\D$ and $k\in\N\cup\{0\}$ we also define spaces $C^{k,\alpha}(G)$ and $C_{\R}^{k,\alpha}(G)$
of $k$ times continuously differentiable functions on $G$, whose all $k$-th derivatives belong to space $C^{\alpha}(G)$ or space $C_{\R}^{\alpha}(G)$.

We also need some algebras of holomorphic functions on $\D$. 
By $A(\D)$ we denote the {\it disc algebra}, that is, the algebra of continuous functions on $\overline{\D}$ which are holomorphic on $\D$,
and by $A^{\alpha}(\D) = A(\D)\cap C^{\alpha}(\overline{\D})$
the algebra of  H\"{o}lder continuous functions on the closed disc which are holomorphic on $\D$.
Using appropriate norms, that is, the maximum norm $\| \cdot\|_{\infty}$ for $A(\D)$ and the H\"{o}lder norm $\|\cdot\|_{\alpha}$ for $A^{\alpha}(\D)$,
these algebras become Banach algebras. Similarly we define $A^{k,\alpha}(\D) = A(\D)\cap C^{k,\alpha}(\overline{\D})$ $(k\in\N\cup\{0\}, 0<\alpha<1)$. 

Recall that Hilbert transform $H$ assigns to a real function $u$ on $\bD$ a real function $Hu$ on $\bD$ such that the harmonic extension
of $f=u+i\, Hu$ to $\D$ is holomorphic on $\D$ and real at $0$.
It is known that $H$ is a bounded linear operator on
$C_{\R}^{k,\alpha}(\bdr\D)$  $(k\in\N\cup\{0\}, 0<\alpha<1)$, \cite[\S 1.6.11]{Weg2}, and hence the harmonic extension 
of $f=u+i\,Hu$ to $\D$ belongs to
$A^{k,\alpha}(\D)$. 
Also, \cite[\S 1.6.11]{Weg2}, the Hilbert transform is a bounded linear operator on the Sobolev space $W^k_p(\bdr\D)$ of $k$ times generalized differentiable functions with derivatives in $L^p(\bdr\D)$  $(k\in\N\cup\{0\}, 1<p<\infty)$ equipped with the norm
\begin{linenomath*}
\begin{equation}
\| f\|_{W_p^k} = \left(\sum_{j=0}^{k} \| D^jf\|_p  \right)^{\frac{1}{p}}.
\end{equation}
\end{linenomath*}

Recall, \cite[\S 1.6.14]{Weg2}, that if $\bdr\D = T_1\cup T_2$ is a partition of $\bdr\D$ in two subarcs $T_1$ and $T_2$ and if $T_0\subseteq T_1$ is a compactly contained subarc of $T_1$, then for $k\in\N\cup\{0\}$,
$1<p<\infty$, $0<\alpha<1$ there exists a constant $C=C(k,p,\alpha)$ such that
\begin{linenomath*}
\begin{equation}\label{est1}
\|Hu\|_{W^k_p(T_0)}\le C (\|u\|_{W^k_p(T_1)} + \|u\|_{L^1(T_2)})
\end{equation}
\end{linenomath*}
and
\begin{linenomath*}
\begin{equation}\label{est2}
\|Hu\|_{C^{k,\alpha}(T_0)}\le C (\|u\|_{C^{k,\alpha}(T_1)} + \|u\|_{L^1(T_2)}).
\end{equation}
\end{linenomath*}
We will also need compact embedding result, \cite[\S 1.1.8]{Weg2}, 
\begin{linenomath*}
\begin{equation}\label{compact}
W^1_p(\bdr\D)\hookrightarrow C^{\beta}(\bdr\D)\hookrightarrow C^{\alpha}(\bdr\D)
\end{equation}
\end{linenomath*}
for $0<\alpha<\beta<1-\frac{1}{p}$, $1<p<\infty$, which holds on arcs in $\bdr\D$ as well.

Since $L\ne\bdr\D$ we can extend $a$ to $\bD$ as a nowhere zero function of class $C^{k+1}$ so
that the winding number $W(a)=0$.
Therefore, \cite[p.\,25]{Weg2}, we can write $\overline{a}$ in the form
\begin{linenomath*}
\begin{equation}
\overline{a} = r e^{h},
\end{equation}
\end{linenomath*}
where $r>0$ is a positive $C^{k,\alpha}$ function on $\bD$ and $h\in A^{k,\alpha}(\D)$.
Thus the original problem  (\ref{eq1L}-\ref{eq1K}) is equivalent to the
problem
\begin{linenomath*}
\begin{equation}\label{eq2L}
{\rm Im}(f_{\ast}(\xi)) = 0 {\rm\ \ for\ \ } \xi\in L
\end{equation}
\end{linenomath*}
and
\begin{linenomath*}
\begin{equation}\label{eq2K}
f_{\ast}(\xi)\in\gamma^{\ast}_{\xi} {\rm\ \ for\ \ } \xi\in\bdr\D\setminus\mathring{L},
\end{equation}
\end{linenomath*}
where $f_{\ast} =  i e^h f$ and $\gamma^{\ast}_{\xi} = i e^{h(\xi)}\gamma_{\xi}$. Observe that the number of zeros of $f_{\ast}$ and $f$
are the same and that $0$ belongs to the interiors of all curves $\gamma^{\ast}_{\xi}$, $\xi\in\bdr\D\setminus\mathring{L}$. Also, since for each $\xi\in\bdr\D$ the transfomation
\begin{linenomath*}
\begin{equation}
w\longmapsto i e^{h(\xi)} w
\end{equation}
\end{linenomath*}
is a composition of a dilation and a rotation, the angle conditions from Theorem \ref{th1} stay the same.

Using a holomorphic automorphism of the unit disc we may even assume that $L=\{ \xi\in\bdr\D; {\rm Im}(\xi)\le 0\}$ is the lower semicircle.
From now on we will consider problem (\ref{eq2L}-\ref{eq2K}) with the addition that $L$ is the lower semicircle and
instead of $f_{\ast}$ and $\zp\gamma^{\ast}_{\xi}\zz_{\xi\in\bD\setminus\mathring{L}}$ we will
still write $f$ and $\zp\gamma_{\xi}\zz_{\xi\in\bD\setminus\mathring{L}}$. 

\begin{remark}\label{rem1}
{\rm One can also create the 'double' of the boundary value problem. 
Using a biholomorphism one can replace the unit disc $\D$ with the upper half-disk $\D_+=\{ \xi\in\D; {\rm Im}(\xi)> 0\}$  and $L$ by the interval $\lbrack -1,1\rbrack$. 

By the reflection principle we see that problem (\ref{eq2L}-\ref{eq2K}) is equivalent to the nonlinear Riemann-Hilbert problem on $\D$, where the boundary curves
$\zp\gamma_{\xi}\zz_{\xi\in\bD_+\setminus\mathring{L}}$ are symmetrically extended and defined on the lower semicircle so that we have
\begin{linenomath*}
\begin{equation}\label{reflected}
\gamma_{\xi} = \overline{\gamma_{\overline{\xi}}}
\end{equation}
\end{linenomath*}
 for every $\xi\in\bdr\D\setminus\{1,-1\}$. In general this symmetrical extension of Jordan curves $\zp\gamma_{\xi}\zz_{\xi\in\bD_+\setminus\mathring{L}}$ to the lower semicircle produces boundary data which are not continuous at points $1$ and $-1$.
Because the biholomorphism from $\D$ to the upper semidisc is in $A^{\frac{1}{2}}(\D)$, we get that the regularity of solutions of 
(\ref{eq2L}-\ref{eq2K}) is in general a half of the regularity of solutions of the symmetrical Riemann-Hilbert problem.}
\end{remark}

We will consider smooth families of smooth Jordan curves  $\zp\gamma_{\xi}\zz_{\xi\in\bdr \D\setminus\mathring{L}}$ in $\C$.
Let $k\in\N$.
The family of Jordan curves $\zp\gamma_{\xi}\zz_{\xi\in\bdr\D\setminus\mathring{L}}$ is a $C^k$ family 
parametrized by $\xi\in\bdr\D\setminus\mathring{L}$ if
there exists a function $\rho\in C^{k}((\bdr\D\setminus\mathring{L})\times\C)$ such that 
\begin{linenomath*}
\begin{equation}\label{eq10}
\gamma_{\xi}=\zp w\in\C; \rho(\xi,w) = 0\zz\ \ \ {\rm and}\ \ \
{\rm Int}(\gamma_{\xi})=\zp w\in\C; \rho(\xi,w) < 0\zz,
\end{equation}
\end{linenomath*}
and the gradient $\frac{\partial\rho}{\partial\overline{w}}(\xi,w)= \rho_{\overline{w}}(\xi,w) \ne 0$ for every $\xi\in\bdr\D\setminus\mathring{L}$ and $w\in\gamma_{\xi}$.
We call $\rho$ a {\it defining function} for $C^k$ family of Jordan curves
$\zp\gamma_{\xi}\zz_{\xi\in\bdr\D\setminus\mathring{L}}$. 
We will consider only bounded families of Jordan curves which all lie in some fixed disc $\overline{\D(0,R)}$, $R>0$, and the 
space $C^{k}((\bdr\D\setminus\mathring{L})\times\overline{\D(0,R)})$ is equipped with the standard $C^k$ norm.

Since we assume that $\gamma_{\pm 1}$ are strongly starshaped Jordan curves, we also assume that for $\rho$, the defining function for Jordan curves $\zp \gamma_{\xi}\zz_{\xi\in\bdr\D\setminus\mathring{L}}$, and $j=\pm 1$ we have
\begin{linenomath*}
\begin{equation}
\rho(j,w) = |w|^2 - R_{j}^2\left(\frac{w}{|w|}\right)
\end{equation}
\end{linenomath*}
for some positive $C^k$ functions $R_j(z)$ on $\C$.

Using parametrization $\theta\mapsto e^{i\theta}$ of the unit circle we will also use the notation $\gamma_{\theta}$, $\rho(\theta,w)$ and $\rho_{\theta}(\theta, w)$
instead of $\gamma_{\xi}$, $\rho(\xi,w)$ and $\rho_{\xi}(\xi,w)$. Also, for a function $h$ on $\bdr\D$, we will write either $h(\xi)$ or $h(\theta)$, where
$\xi=e^{i\theta}$. Observe that if $h$ is holomorphic on $\D$ with well defined derivative on $\bdr\D$, then 
$\frac{\partial h}{\partial\theta} (\theta) = i\xi h'(\xi)$ for $\xi = e^{i\theta}$.

\begin{remark}
{\rm The reflection principle and the symmetric extension to the lower semicircle mentioned in Remark \ref{rem1}
is in terms of defining function $\rho$ given as
\begin{linenomath*}
\begin{equation}
\rho(\xi, w) = \rho(\overline{\xi}, \overline{w})
\end{equation}
\end{linenomath*}
for every $\xi\in\bdr\D\setminus\{1,-1\}$ and every $w\in\C$.}
\end{remark}

\section{Regularity of  solutions}\label{regularity}
In this section we prove regularity of continuous solutions of a specific form of problem (\ref{eq2L}-\ref{eq2K}), where the 
defining function $\rho\in C^{k}((\bdr\D\setminus\mathring{L})\times\C)$ $(k\ge 3)$. 

Let $f\in A(\D)$ be a solution of (\ref{eq2L}-\ref{eq2K}).
It is well known \cite{Chi, Chi-Cou-Suk, For, Weg2} that $f$ restricted to $\bdr\D\setminus\{-1,1\}$ is in $C^{k-1,\alpha}$ for any $0<\alpha<1$. Hence we need information on the regularity of $f$ near points $\xi=\pm 1$. For $j=\pm 1$ we denote $f(j)=w_j\in\R\cap\gamma_j$.

Using M\"{o}bius tranformation from the unit disc $\D$ to the upper half-plane ${\cal H} = \{z\in\C; {\rm Im}(z)>0\}$ we consider the case where $f$ is bounded and continuous on $\overline{\cal H}$ and holomorphic on ${\cal H}$. Also, point $\xi=1$ is mapped into $t=0$ and point $\xi=-1$ into $\infty$.
Now $f$ solves the problem
\begin{linenomath*}
\begin{equation}\label{eq3L}
{\rm Im}(f(t)) = 0 {\rm\ for\ } t\le 0
\end{equation}
\end{linenomath*}
and
\begin{linenomath*}
\begin{equation}\label{eq3K}
f(t)\in\gamma_{t} {\rm\ for\ } t\ge 0.
\end{equation}
\end{linenomath*}
Also, using translation, we will assume that $f(0)=0\in\R\cap\gamma_0$.

Let $\pi \beta_1$ $(\beta_1\in (-1,1)\setminus\{0\})$ be the oriented angle of intersection of the real axis ${\rm Im}(w)=0$ and $\gamma_0$ at $w=f(0)=0$. The orientation of the real axis is positive with respect to the upper half-plane and the orientation of $\gamma_0$ is positive with respect to the interior of $\gamma_0$. Hence $\beta_1> 0$, if the tangent vector to the real axis is rotated counterclockwise by angle $\pi\beta_1$ to get a tangent vector to $\gamma_0$ at point $0$, and $\beta_1<0$, if the tangent vector to the real axis is rotated clockwise by angle $\pi|\beta_1|$ to get a tangent vector to $\gamma_0$ at $0$.

The defining function $\rho$ can near $(0,0)$ for $t\ge 0$ be written as 
\begin{linenomath*}
\begin{equation}\label{rho1}
\rho(t,w) = \rho(0,0)+\rho_t(0,0) t + 2{\rm Re}(\rho_w(0,0) w) + \frac{1}{2}\rho_{tt}(0,0) t^2 +
\end{equation}
\end{linenomath*}
\begin{linenomath*}
\begin{equation}\label{rho2}
+\,  \rho_{w\overline{w}}(0,0)|w|^2 + {\rm Re}(\rho_{ww}(0,0)w^2 + \rho_{tw}(0,0) t w) + 
\sqrt{t^2+|w|^2}\,g(t,w),
\end{equation}
\end{linenomath*}
where $g\in C^{1}(\R\times\C)$ such that $g(0,0) = g_t(0,0)=g_w(0,0)=g_{\overline{w}}(0,0)=0$.

Recall that $\rho(0,0)=0$ and that $\rho_{\overline{w}}(0,0)$ represents an outer normal to $\gamma_0$ at point $w=0$. So we have
\begin{linenomath*}
\begin{equation}
\rho_{\overline{w}}(0,0) = -i \lambda e^{i\pi\beta_1}
\end{equation}
\end{linenomath*}
for some real $\lambda > 0$. We may assume $\lambda = \frac{1}{2}$.

Because
\begin{linenomath*}
\begin{equation}
{\rm Re}(i e^{-i\pi\beta_1} w) = - {\rm Im} (e^{-i\pi\beta_1} w) =  {\rm Im} (e^{i\pi(1-\beta_1)} w)
\end{equation}
\end{linenomath*}
we have 
\begin{linenomath*}
\begin{equation}
\rho(t,w) = A t + {\rm Im}(e^{i\pi(1-\beta_1)} w) +  B t^2 + C |w|^2+
\end{equation}
\end{linenomath*}
\begin{linenomath*}
\begin{equation}
+\,  {\rm Re}(D w^2) + t\,{\rm Re}(E w) +\sqrt{t^2+|w|^2}\,g(t,w)
\end{equation}
\end{linenomath*}
for some $A,B,C\in\R$ and $D,E\in\C$.

Let us assume that we have a solution $f$ of the problem (\ref{eq3L}-\ref{eq3K}) of the form
\begin{linenomath*}
\begin{equation}
f(t) = t^{s} \kappa(t),
\end{equation}
\end{linenomath*}
where $\kappa$ is bounded and continuous on $\overline{\cal H}$, holomorphic on ${\cal H}$, and $0<s<1$ to be determined. 

For $t\le 0$ we have $t=(-1)|t|$ and from (\ref{eq3L}) we get 
\begin{linenomath*}
\begin{equation}
{\rm Im}( e^{i\pi s} \kappa(t))=-{\rm Im}( e^{i\pi (1+s)} \kappa(t))=0.
\end{equation}
\end{linenomath*}
On the other hand for $t>0$ we have
\begin{linenomath*}
\begin{equation}
\frac{1}{t^s}\rho(t,t^s \kappa(s)) = A t^{1-s} + {\rm Im}(e^{i\pi(1-\beta_1)} \kappa(t)) + B t^{2-s} + C t^s|\kappa(t)|^2+
\end{equation}
\end{linenomath*}
\begin{linenomath*}
\begin{equation}
+t^s {\rm Re}(D \kappa(t)^2) + t{\rm Re}(E \kappa(t)) +
\sqrt{t^{2-2s}+|\kappa(t)|^2}\, g(t,t^s\kappa(t))=0.
\end{equation}
\end{linenomath*}
We choose $0<s<1$ so that $\kappa$ solves boundary value problem with continuous boundary data.
That is, we choose $s= 1-\beta_1$, if $\beta_1> 0$, and $s=-\beta_1=|\beta_1|$, if $\beta_1<0$.

Thus $\kappa$ solves the following Riemann-Hilbert problem
\begin{linenomath*}
\begin{equation}\label{eq4L}
{\rm Im}( e^{i\pi (1-\beta_1)} \kappa(t))= 0 {\rm\ for\ } t\le 0
\end{equation}
\end{linenomath*}
and
\begin{linenomath*}
\begin{equation}\label{eq4K}
\widetilde{\rho}(t,\kappa(t)) = 0 {\rm\ for\ } t\ge 0,
\end{equation}
\end{linenomath*}
where, if $\beta_1> 0$,
\begin{linenomath*}
\begin{equation} \label{firsteq1}
\widetilde{\rho}(t,w) =  A t^{\beta_1} + {\rm Im}(e^{i\pi (1-\beta_1)} w) + B t^{1+\beta_1} + C t^{1-\beta_1}|w|^2+
\end{equation}
\end{linenomath*}
\begin{linenomath*}
\begin{equation} \label{firsteq2}
+\,t^{1-\beta_1} {\rm Re}(D w^2) + t{\rm Re}(E w) + 
\sqrt{t^{2\beta_1}+|w|^2}\,g(t,t^{1-\beta_1} w),
\end{equation}
\end{linenomath*}
and, if $\beta_1< 0$,
\begin{linenomath*}
\begin{equation} \label{secondeq1}
\widetilde{\rho}(t,w) =  A t^{1-|\beta_1|} + {\rm Im}(e^{i\pi (1-\beta_1)} w) + B t^{2-|\beta_1|} + C t^{|\beta_1|}|w|^2+
\end{equation}
\end{linenomath*}
\begin{linenomath*}
\begin{equation} \label{secondeq2}
+\,t^{|\beta_1|} {\rm Re}(D w^2) + t{\rm Re}(E w) + 
\sqrt{t^{2-2|\beta_1|}+|w|^2}\,g(t,t^{|\beta_1|} w).
\end{equation}
\end{linenomath*}
For such choice of $s$ are the defining function for problem (\ref{eq4L}-\ref{eq4K}) 
\begin{linenomath*}
\begin{equation}\label{rho-wiggle}
(t,w)\longmapsto
\left\{
\begin{array}{rl}
       \widetilde{\rho}(t,w) = \frac{1}{t^s}\rho(t,t^s w) ;  & t\ge0, w\in\C \\
      {\rm Im}( e^{i\pi(1-\beta_1)} w); & t\le 0, w\in\C
\end{array}
\right.
\end{equation}
\end{linenomath*}
and its partial $w$-derivative
\begin{linenomath*}
\begin{equation}\label{rho-wiggle-w}
(t,w)\longmapsto
\left\{
\begin{array}{rl}
       \widetilde{\rho}_w(t,w) = \rho_w(t,t^s w) ;  & t\ge0, w\in\C \\
      \frac{1}{2i} e^{i\pi(1-\beta_1)}; & t\le 0, w\in\C
\end{array}
\right.
\end{equation}
\end{linenomath*}
continuous on $\R\times\C$. 

On the other hand, the partial derivative of defining function (\ref{rho-wiggle}) with respect to the $t$ variable is not continuous at $t=0$, but, as we will see, it still has certain $L^p$ regularity properties,
which will imply regularity conditions on $\kappa$ and $f$.

We know that $\kappa$ is $C^{k-1,\alpha}$ on $\R\setminus\{0\}$ and we can differentiate (\ref{eq4L}-\ref{eq4K}) on $\R\setminus\{0\}$ to get
\begin{linenomath*}
\begin{equation}\label{eq5L}
{\rm Im}( e^{i\pi  (1-\beta_1)} \kappa'(t))= 0 {\rm\ for\ } t< 0
\end{equation}
\end{linenomath*}
and
\begin{linenomath*}
\begin{equation}\label{eq5K}
\widetilde{\rho}_t (t,\kappa(t))+ 2{\rm Re}(\widetilde{\rho}_w(t,\kappa(t))\kappa'(t)) = 0 {\rm\ for\ } t> 0.
\end{equation}
\end{linenomath*}
For $t>0$ and $\beta_1>0$ we have
\begin{linenomath*}
\begin{equation}\label{pt11}
\widetilde{\rho}_t(t,w)=  A \beta_1 t^{\beta_1-1} + B(1+\beta_1)t^{\beta_1} + (1-\beta_1) C t^{-\beta_1}|w|^2 +
\end{equation}
\end{linenomath*}
\begin{linenomath*}
\begin{equation}\label{pt22}
+(1-\beta_1)t^{-\beta_1} {\rm Re}(D w^2)+ {\rm Re}(E w)+
\frac{\beta_1\, t^{2\beta_1-1}}{\sqrt{t^{2\beta_1}+|w|^2}}\,g(t,t^{1-\beta_1} w) +
\end{equation}
\end{linenomath*}
\begin{linenomath*}
\begin{equation}\label{pt33}
+ \sqrt{t^{2\beta_1}+|w|^2}\, (g_t(t,t^{1-\beta_1}w) + 2{\rm Re}(g_w(t,t^{1-\beta_1}w)(1-\beta_1)t^{-\beta_1}w))
\end{equation}
\end{linenomath*}
and for $t>0$ and $\beta_1<0$ we have
\begin{linenomath*}
\begin{equation}\label{pt111}
\widetilde{\rho}_t(t,w)=  A (1-|\beta_1|) t^{-|\beta_1|} + B(2-|\beta_1|)t^{1-|\beta_1|} + |\beta_1| C t^{|\beta_1|-1}|w|^2
\end{equation}
\end{linenomath*}
\begin{linenomath*}
\begin{equation}\label{pt222}
+|\beta_1|t^{|\beta_1|-1} {\rm Re}(D w^2)+ {\rm Re}(E w)+
\frac{(1-|\beta_1|)\, t^{1-2|\beta_1|}}{\sqrt{t^{2-2|\beta_1|}+|w|^2}}\,g(t,t^{|\beta_1|} w) +
\end{equation}
\end{linenomath*}
\begin{linenomath*}
\begin{equation}\label{pt333}
+ \sqrt{t^{2-2|\beta_1|}+|w|^2}\, (g_t(t,t^{|\beta_1|}w) + 2{\rm Re}(g_w(t,t^{|\beta_1|}w)|\beta_1|t^{|\beta_1|-1}w)).
\end{equation}
\end{linenomath*}
The $t$-derivative of defining function (\ref{rho-wiggle}) is $0$ for $t<0$.

Since $\beta_1\in (-1,1)\setminus\{0\}$ and $\kappa$ is bounded, we have that $\widetilde{\rho}_t(t, \kappa(t))$
is in $L^p_{\rm loc}(\R)$ for 
\begin{linenomath*}
\begin{equation}\label{p-value}
1\le p<\min\left\{\frac{1}{|\beta_1|}, \frac{1}{1-|\beta_1|}\right\} .
\end{equation}
\end{linenomath*}

A similar argument can be used for point $\xi=-1\in\bdr\D$. Let $\pi \beta_{-1}$ $(\beta_{-1}\in (-1,1)\setminus\{0\})$ be the orientied angle of intersection of 
$\gamma_{-1}$ and the real axis  ${\rm Im}(w)=0$ at point $f(-1)$. Now $\beta_{-1}$ is positive, if a positive tangent vector to $\gamma_{-1}$ at $f(-1)$ is rotated counterclockwise to get a positive tangent vector to the real axis and negative otherwise.
For $j=\pm1$ we define $\delta_j = 1-\beta_j$, if $\beta_j\in (0,1)$, and
$\delta_j=|\beta_j|$, if $\beta_j\in (-1,0)$.

To transfer our observations to the boundary value problem (\ref{eq2L}-\ref{eq2K}) on the unit disc,
let $\Psi\in A^{\frac{1}{2}}(\D)$ be a biholomorphic map from $\D$ to the upper half-disc $\D_+$, which maps the lower semicircle $L$ on $\lbrack -1, 1 \rbrack$ so that $\Psi(\pm 1) = \pm 1$.
Let $F(x)=\frac{1}{2}x(3-x^2)$. Then $F(x)-1= -\frac{1}{2}(x-1)^2(x+2)$ and
$F(x)+1=-\frac{1}{2}(x+1)^2(x-2)$. Hence function $\psi(\xi) = F(\Psi(\xi))$ is real on $L$, $\psi(\pm 1) = \pm 1$, and $C^1$ on $\bdr\D$.

Recall that $w_{j}$ is the positive intersection of $\gamma_j$ and the real axis,  $j=\pm 1$. Now we consider only those solutions $f$ of the Cherepanov problem
(\ref{eq2L}-\ref{eq2K}), which are of the form
\begin{linenomath*}
\begin{equation}
f(\xi) = (\xi-1)^{\delta_1}(\xi+1)^{\delta_{-1}} \kappa(\xi) + w_1 \frac{1+\psi(\xi)}{2}+ w_{-1} \frac{1-\psi(\xi)}{2},
\end{equation}
\end{linenomath*}
where $\kappa$ is in $A(\D)$.

We will define two (local) defining functions $\widetilde{\rho}_1(\xi,w)$ for $\xi\ne -1$ and
$\widetilde{\rho}_{-1}(\xi,w)$ for $\xi\ne 1$. Let
\begin{linenomath*}
\begin{equation}
T_1(\xi) = \frac{(\xi-1)}{i(\xi+1)}\ \ \ {\rm and}\ \ \ T_{-1}(\xi) = \frac{1}{T_1(\xi)}=
\frac{i(\xi +1)}{(\xi-1)}.
\end{equation}
\end{linenomath*}
Then $T_1(-i)=T_{-1}(-i)=-1$, and $T_1, T_{-1}$ map the upper semicircle to the positive real axis and the lower semicircle to the negative real axis.
For $j=\pm 1$ and ${\rm Im}(\xi)> 0$ we define
\begin{linenomath*}
\begin{equation}\label{rho1}
\widetilde{\rho}_j(\xi,w)=
\frac{1}{T_j(\xi)^{\delta_j}} \rho\left(\xi,(\xi-1)^{\delta_1} (\xi+1)^{\delta_{-1}}w+ w_1 \frac{1+\psi(\xi)}{2}+ w_{-1} \frac{1-\psi(\xi)}{2}\right)
\end{equation}
\end{linenomath*}
and for ${\rm Im}(\xi)< 0$ we set
\begin{linenomath*}
\begin{equation}\label{rho2}
\widetilde{\rho}_j(\xi,w)=
{\rm Im}\left(e^{i\pi(1-\beta_j)}\frac{(\xi-1)^{\delta_1} (\xi+1)^{\delta_{-1}}}{T_j(\xi)^{\delta_j}}w\right).
\end{equation}
\end{linenomath*}
As before one can check that $\widetilde{\rho}_j$ and $\widetilde{\rho}_{jw}$ are continuous on $\bdr\D\setminus\{-j\}$, $j=\pm1$.
Since $f$ solves the original boundary value problem, we have that $\widetilde{\rho}_j(\xi,\kappa(\xi))=0$, $j=\pm 1$.

Let $\chi:\bdr\D\setminus\{-i\}\rightarrow\lbrack 0,1\rbrack$ be a smooth function such that $\chi(\xi)=1$ for
$\xi=e^{i\theta}$, $-\frac{\pi}{2}<\theta\le\frac{\pi}{3}$ and $\chi(\xi)=0$ for
$\xi=e^{i\theta}$, $\frac{2\pi}{3}\le\theta<\frac{3\pi}{2}$.

We define a new (global) defining function as $\widetilde{\rho}(\xi,w) = \chi(\xi)\widetilde{\rho}_1(\xi,w) + (1-\chi(\xi))\widetilde{\rho}_{-1}(\xi,w)$. Then $\widetilde{\rho}$ and $\widetilde{\rho}_w$ are well defined continuous function on $(\bdr\D\setminus\{-i\})\times\C$. 
If $\beta_1, \beta_{-1}$ have the same sign, then both function are also continuous at $\xi=-i$, but if $\beta_1, \beta_{-1}$
have the opposite signs, then
\begin{linenomath*}
\begin{equation}
\widetilde{\rho}(-i^-,w) = - \widetilde{\rho}(-i^+,w)\ \ {\rm and}\ \ \widetilde{\rho}_w(-i^-,w) = - \widetilde{\rho}_w(-i^+,w),
\end{equation}
\end{linenomath*}
which means that we have a nonorientable bundle as the boundary value data for $\kappa$.

Now locally considered problem (\ref{eq5L}-\ref{eq5K}) for $\kappa(t)$ and $\kappa'(t)$ becomes global boundary value problem for $\kappa(\theta)$ and $\frac{\partial \kappa}{\partial\theta}$ $(\xi=e^{i\theta})$.
Hence $\frac{\partial \kappa}{\partial\theta}$ solves the linear Riemann-Hilbert problem 
\begin{linenomath*}
\begin{equation}
2{\rm Re}\left(\widetilde{\rho}_w(\theta,\kappa(\theta))\frac{\partial \kappa}{\partial\theta}\right)=-\widetilde{\rho}_{\theta}(\theta, \kappa(\theta)),
\end{equation}
\end{linenomath*}
where  $\widetilde{\rho}_w(\theta,\kappa(\theta))$ is either a nonzero continuous function on $\bdr\D$ or 
\begin{linenomath*}
\begin{equation}
\widetilde{\rho}_w(-i^-,\kappa(-i)) = - \widetilde{\rho}_w(-i^+,\kappa(-i))
\end{equation}
\end{linenomath*}
and $\widetilde{\rho}_{\theta}(\theta, \kappa(\theta))$ belongs to the appropriate $L^p(\bdr\D)$ space
\begin{linenomath*}
\begin{equation}\label{allowable}
1\le p < \min\left\{\frac{1}{|\beta_1|}, \frac{1}{1-|\beta_1|}, \frac{1}{|\beta_{-1}|}, \frac{1}{1-|\beta_{-1}|}\right\}.
\end{equation}
\end{linenomath*}
\begin{remark}
{\rm In fact $\widetilde{\rho}_{\theta}(\theta, \kappa(\theta))$ belongs to $L^p_{\rm loc}$ for
\begin{linenomath*}
\begin{equation}
1\le p < \min\left\{\frac{1}{|\beta_1|}, \frac{1}{1-|\beta_1|}\right\}
\end{equation}
\end{linenomath*}
near $\xi=1$ and to $L^p_{\rm loc}$ near $\xi=-1$ for
\begin{linenomath*}
\begin{equation}
1\le p < \min\left\{\frac{1}{|\beta_{-1}|}, \frac{1}{1-|\beta_{-1}|}\right\}.
\end{equation}
\end{linenomath*}
}
\end{remark}

Let $N$ be the winding number of function $\widetilde{\rho}_w(\theta,\kappa(\theta))$, that is, $2N$ is the Maslov index of the associated linear Riemann-Hilbert problem. 
If $\widetilde{\rho}_w(\theta,\kappa(\theta))$ is a continuous function on $\bdr\D$,
Maslov index is an even integer and hence $N$ is an integer. On the other hand, if $\widetilde{\rho}_w(-i^-,\kappa(-i)) = 
- \widetilde{\rho}_w(-i^+,\kappa(-i))$, Maslov index is an odd integer and $N$ is a half of an odd integer.

Let $r(\xi)$ be the square root function, where we take the branch where $\C$ is cut along the negative imaginary axis.
Then function $\widetilde{\rho}_w(\theta,\kappa(\theta))$ can be written in the form
\begin{linenomath*}
\begin{equation}
\widetilde{\rho}_w(\theta,\kappa(\theta)) = \xi^{-N} e^{u+iv}(\theta),
\end{equation}
\end{linenomath*}
where $u$ and $v$ are real continuous functions on $\bdr\D$, \cite[p.\,25]{Weg2}. In the case $N=\frac{2M+1}{2}$, $M\in\Z$, is a half of an odd integer, we define $\xi^N = \xi^M r(\xi)$, which corresponds to the sign changing of $\widetilde{\rho}_w$ at $\xi=-i$.
See also \cite{Cer0, Cer1, Oh}.
Hence $e^{\pm Hv}$ belongs to $L^{p'}(\bdr\D)$ for any $p'\ge 1$, 
\cite[p.\,23]{Weg2} and thus
\begin{linenomath*}
\begin{equation}
e^{\pm i (v+i Hv)}
\end{equation}
\end{linenomath*}
belongs to $L^{p'}(\bdr\D)$ for any $p'\ge 1$. 

Therefore
\begin{linenomath*}
\begin{equation}\label{globalkappa}
{\rm Re}\left(\xi^{-N} e^{ i (v+i Hv)} \frac{\partial\kappa}{\partial\theta}\right) = -e^{-u}e^{-(Hv)}\widetilde{\rho}_{\theta}(\theta, \kappa).
\end{equation}
\end{linenomath*}
We conclude
that the right-hand side belongs to the same $L^p(\bdr\D)$ space as function $\widetilde{\rho}_{\theta}(\theta, \kappa)$.
Since Hilbert transform is bounded in $L^p(\bdr\D)$ spaces, $1<p<\infty$, \cite[p.\,23]{Weg2}, we get that
$\frac{\partial \kappa}{\partial\theta}$ is in $L^p(\bdr\D)$ for the same set  (\ref{allowable}) of values of $p$ as function $\widetilde{\rho}_{\theta}(\theta, \kappa)$.
Therefore $\kappa$ belongs to $L^{1,p}(\bdr\D)$ for all such values of $p$
and this implies that $\kappa\in C^{\beta}(\bdr\D)$, \cite[p.\,10]{Weg2}, where 
\begin{linenomath*}
\begin{equation}\label{alpha1}
0< \beta < \min\{|\beta_1|, 1-|\beta_1|, |\beta_{-1}|, 1-|\beta_{-1}|\}.
\end{equation}
\end{linenomath*}

\begin{remark}
{\rm  Observe that regularity of $\kappa$  and $f$ could also be expressed locally, that is, near $j=\pm 1$  functions $\kappa$ and $f$ belong to H\"{o}lder space $C^{\beta}$, where
$0<\beta<\min\{|\beta_{j}|, 1-|\beta_{j}|\}$.
}
\end{remark}

\begin{proposition}\label{prop1}
Let $k\ge 3$.
Let
$\zp\gamma_{\xi}\zz_{\xi\in\bdr\D\setminus\mathring{L}}$ be a $C^k$ family of Jordan curves in $\C$.
Let $w_j$, $j=\pm 1$, be an intersection of $\gamma_{j}$ and the real axis and let $\pi\beta_j$, $\beta_j\in(-1,1)\setminus\{0\}$, be the oriented angle of intersection of $\gamma_{j}$ with the real axis at point $w_j$. Let
\begin{linenomath*}
\begin{equation}
0< \beta < \min\{|\beta_1|, 1-|\beta_1|, |\beta_{-1}|, 1-|\beta_{-1}|\}.
\end{equation}
\end{linenomath*}
Then for every solution $f$ of  (\ref{eq2L}-\ref{eq2K}) of the form
\begin{linenomath*}
\begin{equation}
f(\xi) = (\xi-1)^{\delta_1}(\xi+1)^{\delta_{-1}} \kappa(\xi) + w_1 \frac{1+\psi(\xi)}{2}+ w_{-1} \frac{1-\psi(\xi)}{2},
\end{equation}
\end{linenomath*}
where $\kappa\in A(\D)$, we have $f,\kappa\in A^{\beta}(\D)$.
\end{proposition}

\begin{remark}
{\rm Observe that in cases where $\beta_1,\beta_{-1}\in (0,1)$, the regularity conditions we get for solutions of the Cherepanov/mixed Riemann-Hilbert problem (\ref{eq1L}-\ref{eq1K}) are consistent with results on the regularity of Riemann maps from the unit disc into simply connected domains bounded by 
Jordan curves which satisfy so called wedge condition, \cite{Les}. If the defining function $\rho$ is independent of $\xi$ and $\beta_j\in(0,1)$, we get $(1-\beta_j)$-regularity. The $\beta_j$-regularity comes from $\xi$-dependence.

Similarly, the expected regularity and the 'order' of zeros of Riemann maps in the cases where $\beta_j\in (-1,0)$ and which are $\xi$ independent, would be $1+|\beta_j|$, but $\xi$-dependence of the defining function $\rho$ changes regularity conditions.

On the other hand, results in \cite{Kai-Leh} show that in the case of nontransversal intersection of the real axis with either $\gamma_1$ or $\gamma_{-1}$ solutions might not be of the form $(\xi-1)^{\delta_1} \kappa(\xi)$ or  $(\xi+1)^{\delta_{-1}} \kappa(\xi)$ for some function $\kappa\in A(\D)$.
}
\end{remark}

\section{Linear Cherepanov boundary value problem}\label{linear-cherepanov}
In this section we consider the linear version of problem (\ref{eq2L}-\ref{eq2K}), that is, a linear Riemann-Hilbert problem with piecewise continuous boundary data, \cite[p.\,169]{Wen}, and $L$ the lower semicircle.
First we consider homogeneous linear problem with piecewise continuous boundary data
\begin{linenomath*}
\begin{equation}\label{eq7L}
{\rm Im}(f(\xi)) = 0 {\rm\ for\ } \xi\in L
\end{equation}
\end{linenomath*}
and
\begin{linenomath*}
\begin{equation}\label{eq7K}
{\rm Re}(\overline{B(\xi)} f(\xi)) = 0 {\rm\ for\ } \xi\in\bdr\D\setminus\mathring{L},
\end{equation}
\end{linenomath*}
where $B$ is a complex nonzero function of class $C^{\beta}$ on the upper semicircle.
The regularity exponent $\beta\in (0,1)$ is bounded by conditions given in Proposition \ref{prop1}.
We may assume without loss of generality that $|B(\xi)|=1$ for all $\xi\in\bdr\D\setminus\mathring{L}$. 

Let $\pi\beta_1$, $\beta_{1}\in (-1, 1)\setminus\{0\}$, be the oriented angle of intersection of the real axis ${\rm Im}(w)=0$ and ${\rm Re}(\overline{B(1)}w) =0$ at point $0$, that is, $B(1)=-i e^{i\pi\beta_1}$. 
Similarly, let $\pi\beta_{-1}$, $\beta_{-1}\in (-1, 1)\setminus\{0\}$, be the oriented angle of intersection of 
${\rm Re}(\overline{B(-1)}w) =0$ and
the real axis ${\rm Im}(w)=0$ at point $0$, that is, $B(-1) =-i  e^{- i\pi\beta_{-1}}$.

We search for solutions $f\in A(\D)$ of (\ref{eq7L}-\ref{eq7K}) of the form $f(\xi)= (\xi-1)^{\delta_1}(\xi+1)^{\delta_{-1}} \kappa(\xi)$
for some $\kappa\in A^{\beta}(\D)$.
Recall that for $j=\pm1$ we defined $\delta_j = 1-\beta_j$, if $\beta_j\in (0,1)$, and
$\delta_j=|\beta_j|$, if $\beta_j\in (-1,0)$. Hence we also have $f\in A^{\beta}(\D)$.

To define noninteger powers of $(\xi-1)$ and $(\xi+1)$ we take appropriate branches of the complex logarithm. For $(\xi-1)^{\delta_1}$ the complex plane is cut along positive real numbers so  that the argument of $(\xi-1)$ for $\xi\in\bdr\D$ lies on interval $(\frac{\pi}{2},\frac{3\pi}{2})$, and for $(\xi+1)^{\delta_{-1}}$ the complex plane is cut along negative real numbers and the argument of $(\xi+1)$  for $\xi\in\bdr\D$ lies on interval $(-\frac{\pi}{2},\frac{\pi}{2})$.

An argument similar to the argument in Section \ref{regularity} shows that $\kappa$ solves homogeneous linear Riemann-Hilbert problem
\begin{linenomath*}
\begin{equation}\label{linear-widetilde}
{\rm Re}(\overline{\widetilde{B}(\xi)}\kappa(\xi))=0\ {\rm for\ all\ }\xi\in\bdr\D,
\end{equation}
\end{linenomath*}
where $\overline{\widetilde{B}}\in C^{\beta}(\bdr D\setminus\{1\})$ is defined as
\begin{linenomath*}
\begin{equation} \label{linear-widetilde1}
\overline{\widetilde{B}(\xi)}=
\left\{
\begin{array}{rl}
       \overline{B(\xi)}\left(\frac{\xi-1}{|\xi-1|}\right)^{\delta_1} \left(\frac{\xi+1}{|\xi+1|}\right)^{\delta_{-1}}, & {\rm if\  } {\rm Im}(\xi) > 0 \\
     \pm i\left(\frac{\xi-1}{|\xi-1|}\right)^{\delta_1}\left(\frac{\xi+1}{|\xi+1|}\right)^{\delta_{-1}}, & {\rm if\ } {\rm Im}(\xi)<0 
\end{array}
\right.
\end{equation}
\end{linenomath*}
with the left and the right limits at $\xi=\pm1$. The sign for ${\rm Im}(\xi)<0$ is chosen so that $\widetilde{B}$ is continuous at $-1$, that is, we have plus sign, if $\beta_{-1}<0$,
and minus sign, if $\beta_{-1}>0$.
At point $\xi=1$ function $B$ might not be continuous.  In general we have
$\widetilde{B}(1^+)=\pm \widetilde{B}(1^-)$.
See  \cite[p.\,169-170]{Wen} for more.

Each factor 
\begin{linenomath*}
\begin{equation}
\frac{\xi-1}{|\xi-1|}, \frac{\xi+1}{|\xi+1|}
\end{equation}
\end{linenomath*}
changes the argument by $\pi$ when $\xi$ passes $\bdr\D$ once in the positive direction.
Hence possible widing number of $\widetilde{B}$ is either an integer (Maslov index of problem (\ref{linear-widetilde}) is even) or
a half of an odd integer  (Maslov index of problem (\ref{linear-widetilde}) is odd).

\begin{example}
{\rm Consider the case $B(e^{i\theta})= e^{i \theta}$ for $\theta\in [0,\pi]$. In particular we have $\beta_1=\beta_{-1}=\frac{1}{2}$ . Then we get
\begin{linenomath*}
\begin{equation}
\widetilde{B}(\xi)=
\left\{
\begin{array}{rl}
        e^{-i\frac{\pi}{4}}B(\xi)\,\overline{\xi}^{\frac{1}{2}}= e^{i\frac{2\theta-\pi}{4}},  & {\rm if\  } 0\le \theta\le\pi \\
     e^{i\frac{3\pi}{4}}\,\overline{\xi}^{\frac{1}{2}} =e^{i\frac{3\pi-2\theta}{4}} , & {\rm if\ }  \pi< \theta<2\pi.
\end{array}
\right.
\end{equation}
\end{linenomath*}
Hence the winding number $W(\widetilde{B})=0$.
Using identification of the boundary problem (\ref{eq7L}-\ref{eq7K}) with the problem on the unit disc with reflected boundary conditions (\ref{reflected}), this example corresponds to the linearization of the boundary value problem, where all boundary curves are unit circles and we linearize at $f(z)=z$. The family of (nearby) solutions which are real on the real axis is one-dimensional
$f_a(z) = \frac{z-a}{1-az}$, where $a\in (-1,1)$ is a real number.
}
\end{example}

\begin{example}
{\rm Consider the case $B(e^{i\theta})= 1$ for $\theta\in [0,\pi]$. In particular we have $\beta_1=-\beta_{-1}=\frac{1}{2}$ . Then we get
\begin{linenomath*}
\begin{equation}
\widetilde{B}(\xi)=  e^{-i\frac{\pi+2\theta}{4}}
\end{equation}
and the winding number $W(\widetilde{B})=-\frac{1}{2}$.
Using identification of the boundary problem (\ref{eq7L}-\ref{eq7K}) with the problem on the unit disc with reflected boundary conditions (\ref{reflected}), this example corresponds to the linearization of the 
problem where all boundary curves are unit circles and we linearize at function $f(z)=1$. The family of (nearby) solutions which are real on the real axis is zero-dimensional.
\end{linenomath*}
 }
\end{example}

The dimension of the space of solutions in $A^{\beta}(\D)$ depends on the winding number $W(\widetilde{B})$ of function $\widetilde{B}$. It equals $2 W(\widetilde{B})+1$ if $W(\widetilde{B})\ge -\frac{1}{2}$, see \cite[p.\,25, p.\,59]{Weg2} and \cite{Cer0, Cer1, Oh}.
We define the winding number of $B\in C^{\alpha}(\bdr\D\setminus\mathring{L})$ as the winding number of
$\widetilde{B}$.

Now we can solve appropriate nonhomogeneous linear Riemann-Hilbert problem with piecewise continuous boundary data
\begin{linenomath*}
\begin{equation}\label{eq8L}
{\rm Im}(f(\xi)) = 0 {\rm\ for\ } \xi\in L
\end{equation}
\end{linenomath*}
and
\begin{linenomath*}
\begin{equation}\label{eq8K}
{\rm Re}(\overline{B(\xi)} f(\xi)) = b(\xi) {\rm\ for\ } \xi\in\bdr\D\setminus\mathring{L},
\end{equation}
\end{linenomath*}
where $B$ is as above and $b$ a real function on $\bdr\D$ of the form
\begin{linenomath*}
\begin{equation}
b(\xi) = |\xi-1|^{\delta_1}|\xi+1|^{\delta_{-1}} \widetilde{b}(\xi)
\end{equation}
\end{linenomath*}
for some function $\widetilde{b}\in C_{\R}^{\beta}(\bdr\D)$ which equals $0$ on $L$.

To solve (\ref{eq8L}-\ref{eq8K}) in the space of functions $f\in A^{\beta}(D)$ of the form $f(\xi)= (\xi-1)^{\delta_1}(\xi+1)^{\delta_{-1}} \kappa(\xi)$ for some $\kappa\in A^{\beta}(\D)$ is equivalent to solve the problem
\begin{linenomath*}
\begin{equation}\label{wiggle-B}
{\rm Re}(\overline{\widetilde{B}}(\xi)\kappa(\xi))=\widetilde{b}(\xi)\ {\rm for\ all\ }\xi\in\bdr\D.
\end{equation}
\end{linenomath*}
It is well known that if $W(B)=W(\widetilde{B})\ge -\frac{1}{2}$, then the equation is solvable for any $\widetilde{b}\in C_{\R}^{\beta}(\bdr\D)$,  see \cite[p.\,25, p.\,59]{Weg2} and \cite{Cer0, Cer1, Oh}.

\begin{remark}
{\rm If the winding number $W(B)=W(\widetilde{B})$ is an odd integer, the function on the right-hand side of (\ref{wiggle-B}) needs to 
belong to a special space of H\"{o}lder continuous real functions on $\bdr\D\setminus\{1\}$ of the form 
$b_0(r(\xi))$, where $r(\xi)$ is the principal branch of the square root and $b_0\in C^{\beta}_{\R}(\bdr\D)$ is an odd function. 
Hence we need condition $\widetilde{b}(1^-) + \widetilde{b}(1^+) =0$, which is satisfied because
in our case we have $\widetilde{b}(1^-) = \widetilde{b}(1^+) =0$.
See \cite{Cer0, Cer1} for more information.}
\end{remark}

\begin{proposition}\label{prop2}
Let $0<\beta<1$. Let $B : \bdr\D\setminus\mathring{L}\rightarrow\C\setminus\{0\}$
be a non-vanishing complex function in $C^{\beta}(\bdr\D\setminus\mathring{L})$ and let $W(B)\ge -\frac{1}{2}$.
Then for every real function $b$ on $\bdr\D$ of the form
\begin{linenomath*}
\begin{equation}
b(\xi) = |\xi-1|^{\delta_1}|\xi+1|^{\delta_{-1}} \widetilde{b}(\xi)
\end{equation}
\end{linenomath*}
for some $\widetilde{b}\in C^{\beta}_{\R}(\bdr\D)$ which equals $0$ on $L$,
there exists a solution $f$ of the linear Cherepanov problem
\begin{linenomath*}
\begin{equation}\label{eq9L}
{\rm Im}(f(\xi)) = 0 {\rm\ for\ } \xi\in L
\end{equation}
\end{linenomath*}
and
\begin{linenomath*}
\begin{equation}\label{eq9K}
{\rm Re}(\overline{B(\xi)} f(\xi)) = b(\xi) {\rm\ for\ } \xi\in\bdr\D\setminus\mathring{L}
\end{equation}
\end{linenomath*}
 of the form
\begin{linenomath*}
\begin{equation}
f(\xi) = (\xi-1)^{\delta_1}(\xi+1)^{\delta_{-1}} \kappa(\xi),
\end{equation}
\end{linenomath*}
where $\kappa\in A^{\beta}(\D)$. Moreover, the space of solutions of this form is $2W(B)+1$ dimensional real subspace of $A^{\beta}(\D)$.
\end{proposition}

\begin{proposition}\label{prop3}
Let $\zp\gamma_{\xi}\zz_{\xi\in\bD\setminus\mathring{L}}$ be a $C^k$ $(k\ge 3)$ family of Jordan curves in $\C$ and let $\rho_0\in C^k((\bdr\D\setminus\mathring{L})\times\C)$ be its defining function.
Let $\beta_1,\beta_{-1}$ and $\beta$ be as in Proposition \ref{prop1}.
Let $f_0$ be a solution of the Cherepanov problem  (\ref{eq2L}), (\ref{eq2K}) of the form
\begin{linenomath*}
\begin{equation}
f_0(\xi) = (\xi-1)^{\delta_1}(\xi+1)^{\delta_{-1}} \kappa_0(\xi) + w_1 \frac{1+\psi(\xi)}{2}+ w_{-1} \frac{1-\psi(\xi)}{2},
\end{equation}
\end{linenomath*}
where $\kappa_0\in A^{\beta}(\D)$.
Then the mapping $\Phi(\kappa) : A^{\beta}(\D)\rightarrow C^{\beta}_{\R}(\bdr\D)$, for each $\kappa$ evaluated at point $\xi\in\bdr\D$ as
\begin{linenomath*}
\begin{equation}\label{map-1}
\left\{
\begin{array}{ll}
     \rho_0(\xi, (\xi-1)^{\delta_1}(\xi+1)^{\delta_{-1}} \kappa(\xi) + w_1 \frac{1+\psi(\xi)}{2}+ w_{-1} \frac{1-\psi(\xi)}{2}), & {\rm if\  } {\rm Im}(\xi) \ge 0 ,\\
   \pm {\rm Im}((\xi-1)^{\delta_1}(\xi+1)^{\delta_{-1}} \kappa(\xi) + w_1 \frac{1+\psi(\xi)}{2}+ w_{-1} \frac{1-\psi(\xi)}{2}), & {\rm if\ } {\rm Im}(\xi) < 0 
\end{array}
\right.
\end{equation}
\end{linenomath*}
is differentiable at $\kappa_0$ with the derivative $(D\Phi)(\kappa_0)$ acting on $\kappa\in A^{\beta}(\D)$ as
\begin{linenomath*}
\begin{equation}\label{map-2}
\left\{
\begin{array}{ll}
     2{\rm Re}(\partial\rho_{0w}(\xi, f_0(\xi)) (\xi-1)^{\delta_1}(\xi+1)^{\delta_{-1}} \kappa(\xi)), & {\rm if\  } {\rm Im}(\xi) \ge 0,\\
   \pm {\rm Im}((\xi-1)^{\delta_1}(\xi+1)^{\delta_{-1}} \kappa(\xi)), & {\rm if\ } {\rm Im}(\xi) < 0.
\end{array}
\right.
\end{equation}
\end{linenomath*}
The sign for ${\rm Im}(\xi)<0$ is chosen as in (\ref{linear-widetilde}-\ref{linear-widetilde1}).
\end{proposition}

\begin{remark}
{\rm Let $\Omega\subset C^{k+1} ((\bdr\D\setminus{\mathring{L}})\times\C)$ be an open subset of defining functions $\rho$ of the families of Jordan curves over $\bdr\D\setminus{\mathring{L}}$ such that
the intersection of the corresponding $\gamma_1$ and $\gamma_{-1}$ with the real axis at some points $w_1\in\gamma_1$ and $w_{-1}\in\gamma_{-1}$ are transversal with the oriented angles of intersection given by
$\beta_1, \beta_{-1}\in(-1,1)\setminus\{0\}$. Then, at least locally, $w_1, w_{-1}$ and $\beta_1,\beta_{-1}$ smoothly depend on $\rho$.
Let 
\begin{linenomath*}
\begin{equation}
X =\{ (\kappa,\rho)\in A^{\beta}(\D)\times\Omega; {\rm Im}((\xi-1)^{\delta_1}(\xi+1)^{\delta_{-1}} \kappa(\xi)) = 0, {\ \rm if\ } {\rm Im}(\xi) < 0\}
\end{equation}
\end{linenomath*}
which is a Banach submanifold of $A^{\beta}(\D)\times\Omega$. 
Also, let 
\begin{linenomath*}
\begin{equation}
Y =\{b(\xi) = |\xi-1|^{\delta_1}|\xi+1|^{\delta_{-1}} \widetilde{b}(\xi);  \widetilde{b}\in C^{\beta}_{\R}(\D),  \widetilde{b}(\xi)=0, {\ \rm if\ } {\rm Im}(\xi) < 0\}.
\end{equation}
\end{linenomath*}
The mapping $\Phi : X\rightarrow Y$ defined as in (\ref{map-1}) has partial derivative with respect to $\kappa$ as a map from
\begin{linenomath*}
\begin{equation}
X_{\rho} =\{ \kappa\in A^{\beta}(\D); {\rm Im}((\xi-1)^{\delta_1}(\xi+1)^{\delta_{-1}} \kappa(\xi)) = 0, {\ \rm if\ } {\rm Im}(\xi) < 0\}
\end{equation}
\end{linenomath*}
to $Y$ of the form (\ref{map-2}). If the winding number $W(B)$ of the Cherepanov problem defined by (\ref{map-2}) is greater or equal to $-\frac{1}{2}$, then the partial derivative is surjective with
$2W(B)+1$ dimensional kernel.
Hence implicit function theorem applies and there is a neighbourhood of $\rho_0$ in $\Omega$  and a neighbourhood of $\kappa_0$ in $A^{\beta}(\D)$
such that for every $\rho\in\Omega$ close to $\rho_0$ there is a $2W(B)+1$ dimensional family of solutions of (\ref{eq2L}-\ref{eq2K}) near $\kappa_0$.
}
\end{remark}

\section{A priori estimates} \label{apriori}
\subsection{A priori estimates on function $f$}
To get existence results using continuity method we need a priori estimates on solutions of (\ref{eq1L}-\ref{eq1K}). It is well known that such a priori estimates can only be achieved for the family of solutions with no zeros on $\D$, \cite{For, Sni, Weg2}. We follow the approach in \cite{For}.

By assumption all Jordan curves $\{\gamma_{\xi}\}_{\xi\in\bdr\D\setminus\mathring{L}}$ contain point $0$ in their interiors. 
Hence the function
\begin{linenomath*}
\begin{equation}
(\theta, w)\longmapsto w\rho_w(\theta,w), 
\end{equation}
\end{linenomath*}
defined for $(\theta, w)$ such that  $w\in\gamma_{\theta}$,
is homotopic to $0$ in $\C\setminus\{0\}$ and it can be written in the form
\begin{linenomath*}
\begin{equation}\label{expform}
w\rho_w(\theta,w) = e^{c(\theta,w) +  i  d(\theta,w) }
\end{equation}
\end{linenomath*}
for some $C^{k-1}$ functions $c$ and $d$, defined for $(\theta, w)$ such that $\theta\in\lbrack 0, \pi\rbrack$ and  $w\in\gamma_{\theta}$.
Observe that for each $w\in\gamma_{\theta}$ function $d(\theta,w)$ represents the angle between $w$ and the outer normal to $\gamma_{\theta}$ at point $w$.
\begin{remark}
{\rm There exists a $C^k$ isotopy $\rho^t$, $t\in\lbrack 0,1\rbrack$, where $\rho^0=\rho$ and $\rho^1 (\xi,w) = |w|^2 - R^2$  for $R>0$ large enough, such that the gradient $\rho^t_{\overline{w}}$ is nonzero on $\rho^t=0$ for each $t$, \cite{For}. Then one can find 
$C^{k-1}$ functions $c(t,\theta,w)$ and $d(t,\theta,w)$ such that (\ref{expform}) holds for each $t\in\lbrack 0,1\rbrack$ and
$w\in\gamma^t_{\theta}$. In addition, the isotopy can be made such that for every $t\in\lbrack 0, 1\rbrack$, $j=\pm 1$,
Jordan curves $\gamma^t_{j}$ are strongly starshaped with respect to $0$ and that for each 
$w\in \gamma^t_{\omega_{j}}$ the angle between $w$ and the normal to $\gamma^t_{\omega_{j}}$
at $w$ is less than $\frac{\pi}{10}$.
}
\end{remark}

Instead of solving  (\ref{eq2L}-\ref{eq2K}) on the unit disc we consider equivalent problem on the upper semidisc $\D^+ =\{ z\in\D; {\rm Im}(z)>0\}$, where the role of the lower semicircle $L$ is replaced by the interval $\lbrack -1,1\rbrack$.
Using the reflection principle  $f(\xi) = \overline{f(\overline{\xi})})$ we can holomorphically extend every solution $f$ of  (\ref{eq2L}-\ref{eq2K}) to the unit disc such that it solves nonlinear Riemann-Hilbert problem defined by the function $\rho$ which we get as an extension of the original function $\rho$ using the reflection to the lower semicircle as
\begin{linenomath*}
\begin{equation}
\rho(\xi, w) = \rho(\overline{\xi},\overline{w}) {\rm\ for\ } \xi\ne \pm 1.
\end{equation}
\end{linenomath*}
For $\xi=\pm 1$ function $\rho(\xi, w)$ has well defined limits as $\xi$ approaches $\pm 1$ from above and below.
Then we have
\begin{linenomath*}
\begin{equation}\label{symmetry}
\rho_w(\xi, w) = \rho_{\overline{w}}(\overline{\xi},\overline{w}) = \overline{\rho_w(\overline{\xi},\overline{w})} {\rm\ for\ } \xi\ne \pm 1
\end{equation}
\end{linenomath*}
and hence 
\begin{linenomath*}
\begin{equation}
\rho_w(\xi, w) w = \rho_{\overline{w}}(\overline{\xi},\overline{w}) w= \overline{\rho_w(\overline{\xi},\overline{w}) \overline{w}} .
\end{equation}
\end{linenomath*}
Therefore  $c(\overline{\xi},\overline{w}) = c(\xi,w)$ and $d(\overline{\xi},\overline{w}) = - d(\xi,w)$.
Also, observe that for $w$, an intersection of $\gamma_1$ with the real axis, we have 
\begin{linenomath*}
\begin{equation}
\rho_w(1+, w) = \overline{\rho_w(1-,\overline{w})} =  \overline{\rho_w(1-,w)}
\end{equation}
\end{linenomath*}
and similarly for an intersection of $\gamma_{-1}$ with the real axis.

Thus for every solution $f$ of (\ref{eq2L}-\ref{eq2K}) the absolute value of $f(\theta)\rho_w(\theta, f(\theta))$ is well defined and continuous on $\bdr\D$, whereas 
\begin{linenomath*}
\begin{equation}
d(0+,f(0+))= - d(2\pi-, f(2\pi-))
\end{equation}
\end{linenomath*}
and similarly at $\theta=\pi$.

Let $f$ be a solution of the symmetrized boundary value problem
with no zeros. Hence $f$ can be written in the exponential form
\begin{linenomath*}
\begin{equation}
f = e^g.
\end{equation}
\end{linenomath*}
\begin{remark}\label{remark52}
{\rm Since the biholomorphic map $\psi$ from $\D$ to the upper half-disc $\D_+$ is of class $C^{\frac{1}{2}}$,
a $C^{\beta}$ estimate on solutions of the symmetrized boundary value problem gives $C^{\frac{\beta}{2}}$ estimate on solutions of
(\ref{eq2L}-\ref{eq2K}).}
\end{remark}
Let us differentiate function $\rho(\theta, f(\theta))$
to get
\begin{linenomath*}
\begin{equation}
\rho_{\theta}(\theta, f(\theta)) + 2 {\rm Re}(\rho_w(\theta, f(\theta)) \frac{\partial f}{\partial\theta}(\theta)) = 0.
\end{equation}
\end{linenomath*}
Since $f=e^g$, we get
\begin{linenomath*}
\begin{equation}
\rho_{\theta}(\theta, f(\theta)) + 2 {\rm Re}(\rho_w(\theta, f(\theta)) f(\theta) \frac{\partial g}{\partial\theta}(\theta)) = 0
\end{equation}
\end{linenomath*}
and so
\begin{linenomath*}
\begin{equation}
\rho_{\theta}(\theta, f(\theta)) + 2 {\rm Re}(e^{c(\theta,f(\theta))+ i d(\theta,f(\theta))} \frac{\partial g}{\partial\theta}(\theta)) = 0.
\end{equation}
\end{linenomath*}
From here we get
\begin{linenomath*}
\begin{equation} \label{equation1}
2 {\rm Re}(e^{i (d(\theta,f(\theta))+ i Hd(\theta,f(\theta)))} \frac{\partial g}{\partial\theta}(\theta)) = -\rho_{\theta}(\theta, f(\theta)) e^{-c(\theta,f(\theta))-Hd(\theta,f(\theta))}.
\end{equation}
\end{linenomath*}
Observe that function
\begin{linenomath*}
\begin{equation} \label{holfun}
\theta\longmapsto e^{i (d(\theta,f(\theta))+ i Hd(\theta,f(\theta)))} \frac{\partial g}{\partial\theta}(\theta)
\end{equation}
\end{linenomath*}
extends holomorphically to the unit disc with value $0$ at $0$. 

We will get $C^{\beta}$ a priori estimates on $g$ and hence on $f$ by getting $C^{\beta}$ a priori estimates on function (\ref{holfun}).
Using Hilbert transform it is enough to get $C^{\beta}$  a priori estimates on its real part. Hence we need a priori estimates on the right hand side of (\ref{equation1}).

Function
\begin{linenomath*}
\begin{equation}
\theta\longmapsto  -\rho_{\theta}(\theta, f(\theta)) e^{-c(\theta,f(\theta))}
\end{equation}
\end{linenomath*}
is bounded with the bound which does not depend on function $f$ but only on the data $\gamma_{\xi\in\bdr\D\setminus\mathring{L}}$ and defining function $\rho$. The bound can also be found to be independent of the $C^k$ isotopy $\rho^t$, $t\in\lbrack 0, 1\rbrack$.
Hence one needs a priori bound on function 
\begin{linenomath*}
\begin{equation} \label{ehd}
\theta\longmapsto  e^{\pm Hd(\theta,f(\theta))}.
\end{equation}
\end{linenomath*}
Recall that, \cite[p.\,23]{Weg2}, for $u\in L^{\infty}(\bdr\D)$, such that $\| u\|_{\infty}<\frac{\pi}{2p}$ $(1\le p<\infty)$ we have the estimate
\begin{linenomath*}
\begin{equation}\label{zygmund}
\|e^{Hu}\|_p\le\left (\frac{2\pi}{\cos(p\|u\|_{\infty})}\right)^{\frac{1}{p}}.
\end{equation}
\end{linenomath*}

Let $a\in(0,\frac{\pi}{5})$ and let $\chi_0, \chi_{\pi}$ be smooth functions on $\lbrack 0,\pi\rbrack$ with values in $\lbrack 0,1\rbrack$ such that $\chi_0(t)=1$ on $\lbrack 0,a\rbrack$,  $\chi_{\pi}(t)=1$ on $\lbrack \pi-a, \pi\rbrack$,  $\chi_0(t)=0$ on  $\lbrack 2a, \pi\rbrack$, and $\chi_{\pi}(t)=0$ on $\lbrack 0, \pi-2a\rbrack$.

Let us consider the function
\begin{linenomath*}
\begin{equation} \label{contd}
\widetilde{d}(\theta,w) = d(\theta,w) - \chi_0(\theta)d_0(w) - \chi_{\pi}(\theta)d_{\pi}(w)
\end{equation}
\end{linenomath*}
for $\theta\in\lbrack 0,\pi\rbrack$ and
$\widetilde{d}(\theta,w) =-\widetilde{d}(2\pi-\theta,\overline{w})$ for $\theta\in\lbrack \pi,2\pi\rbrack$.
Here we used notation $d_0(w)=d(0+,w)$ and $d_{\pi}(w)=d(\pi-,w)$.

We see that $\widetilde{d}(0,w) = \widetilde{d}(\pi,w)=0$ and so $\widetilde{d}(\theta,w)$ is a continuous function on $\bdr\D\times\C$.
Let $1< \widetilde{p}<\infty$ be given. By results from \cite[p.\,881]{For} we can write $\widetilde{d} = {\rm Re}(q) + \widetilde{e}$, where $ \widetilde{p}\|\widetilde{e}\|_{\infty}<\frac{\pi}{2}$ and
$q$ is a finite sum of terms of the form $e^{ij\theta} w^m$, $j\in\Z$, $m\in\N\cup\{0\}$, on which Hilbert transform acts as a bounded nonlinear operator from $A(\bdr\D)$ into $C(\bdr\D)$.

Therefore for a given solution $f$ of (\ref{eq2L}-\ref{eq2K}) with no zeros we can write continuous function $\widetilde{d}(\theta,f(\theta))$ on $\bdr\D$ 
in the form
\begin{linenomath*}
\begin{equation}
\widetilde{d}(\theta,f(\theta)) = {\rm Re}(q(\theta,f(\theta))) + \widetilde{e}(\theta,f(\theta))
\end{equation}
\end{linenomath*}
and so 
\begin{linenomath*}
\begin{equation}
H(\widetilde{d}) = H({\rm Re}(q)) + H(\widetilde{e})
\end{equation}
\end{linenomath*}
where the first term is uniformly bounded and for the second we have 
$\|\widetilde{e}\|_{\infty}<\frac{\pi}{2\widetilde{p}}$.
Hence
\begin{linenomath*}
\begin{equation}\label{2and3}
e^{\pm H\widetilde{d}} = e^{\pm H{\rm Re}(q)} e^{\pm H\widetilde{e}}
\end{equation}
\end{linenomath*}
where the first factor is uniformly bounded and the second factor is bounded in $L^{ \widetilde{p}}(\bdr\D)$ for a given $1< \widetilde{p}<\infty$. 

Since for a given $\widetilde{p}\in (1,\infty)$ we can get $L^{\widetilde{p}}(\bdr\D)$ bounds on (\ref{2and3}), the boundedness of $e^{\pm Hd}$ in some $L^{p}(\bdr\D)$ is determined by Hilbert transform of the extension of function
$\chi_0(\theta)d_0(f(\theta))+\chi_{\pi}(\theta)d_{\pi}(f(\theta))$  to $\lbrack 0,2\pi\rbrack$.

Recall that $\gamma_{\pm 1}$ are strongly starshaped Jordan curves with respect to $0$
and we may assume that for $j=\pm 1$ we have
\begin{linenomath*}
\begin{equation}
\rho(j,w) = |w|^2 - R_{j}^2\left(\frac{w}{|w|}\right)
\end{equation}
\end{linenomath*}
for some positive $C^k$ function $R_j(z)$ on $\C$.
A short calculation gives
\begin{linenomath*}
\begin{equation}
\rho_{w}(j,w)=\overline{w}-2 R_j \left(-\frac{1}{2}\frac{\overline{w}^2}{|w|^3} (R_j)_{\overline{z}} + \frac{1}{2}\frac{1}{|w|} (R_j)_{z}\right)
\end{equation}
\end{linenomath*}
and so
\begin{linenomath*}
\begin{equation}\label{starshaped}
w\rho_{w}(j,w)=|w|^2- 2 i\frac{R_j}{|w|} {\rm Im}( w(R_j)_{z})
\end{equation}
\end{linenomath*}
which has strictly positive real part on 
$\zp\gamma_{\xi}\zz_{\xi\in\bD\setminus\mathring{L}}$.
Functions $d_0$ and $d_{\pi}$ represent the argument of (\ref{starshaped}). 
By compactness it follows that there exists $0<\beta_0<1$, such that
$|d_0(w)|\le \frac{\pi}{2}\beta_0$ and $|d_{\pi}(w)|\le \frac{\pi}{2}\beta_0$ on 
$\zp\gamma_{\xi}\zz_{\xi\in\bD\setminus\mathring{L}}$ and therefore 
\begin{linenomath*}
\begin{equation}
|\chi_0(\theta)d_0(f(\theta)) + \chi_{\pi}(\theta)d_{\pi}(f(\theta))|\le  \frac{\pi}{2}\beta_0
\end{equation}
\end{linenomath*}
for every $\theta$. 
Also, if there is an open condition on the size of $d_{j\pi}$ on $\gamma_j$, such as
$|d_{i\pi}(w)|<\frac{\pi}{2}\beta_0$, $j=\pm 1$, we can, by choosing the supports of functions $\chi_0$ and $\chi_{\pi}$ small enough, that is, by choosing $a>0$ small enough, assume that the same condition on the size holds for function
$\chi_0(\theta)d_0(f(\theta)) + \chi_{\pi}(\theta)d_{\pi}(f(\theta))$ for all $\theta$.
Observe also that $|d_0(f(0))|= \pi|\beta_1-\frac{1}{2}|$ and $|d_{\pi}(f(\pi))|=\pi||\beta_{-1}|-\frac{1}{2}|$.
and so
\begin{linenomath*}
\begin{equation}
\left||\beta_j|-\frac{1}{2}\right|<\frac{\beta_0}{2}, \ \ \ j=\pm 1.
\end{equation}
\end{linenomath*}

By (\ref{zygmund}) and (\ref{2and3}) we get that for every fixed $1<p<\infty$ such that $p\beta_0<1$ the estimate
\begin{linenomath*}
\begin{equation}
\|e^{\pm Hd}\|_p\le C
\end{equation}
\end{linenomath*}
holds. Hence we also have a priori $L^p$ estimate on function (\ref{holfun}).

Since 
\begin{linenomath*}
\begin{equation}
\frac{\partial f}{\partial\theta} =f \frac{\partial g}{\partial\theta},
\end{equation}
\end{linenomath*}
an estimate on $\frac{\partial g}{\partial\theta}$ will give an estimate on $\frac{\partial f}{\partial\theta}$.
We can write
\begin{linenomath*}
\begin{equation}
\frac{\partial g}{\partial\theta} = \left(e^{-i (d + i Hd)}\right)\left(e^{i (d+ i Hd)} \frac{\partial g}{\partial\theta}\right).
\end{equation}
\end{linenomath*}
By assumptions of Theorem \ref{th1} we have $\beta_0\le\frac{1}{5}<\frac{1}{2}$ and we can choose $p>2$. By Cauchy-Schwarz inequality we then have
\begin{linenomath*}
\begin{equation}
\left\|\frac{\partial g}{\partial\theta}\right\|_{\frac{p}{2}}\le \left\|e^{-i (d + i Hd)}\right\|_p\left\|e^{i (d+ i Hd)} \frac{\partial g}{\partial\theta}\right\|_p.
\end{equation}
\end{linenomath*}
From here we get $L^{\frac{p}{2}}$ a priori estimates on $\frac{\partial g}{\partial\theta}$ which imply a priori estimates on $g$ and $f$ in H\"{o}lder space $C^{\beta}(\bdr\D)$ for $0<\beta<1-\frac{2}{p}<2(\frac{1}{2}-\beta_0)$.
Recall (Remark \ref{remark52}) that this gives H\"{o}lder space a priori estimates on solutions  with no zeros of the nonsymmetrical problem (\ref{eq2L}-\ref{eq2K}) for $\beta\in(0,\frac{1}{2}-\beta_0)$.

\subsection{A priori estimates on function $\kappa$}

We also need a priori estimates on function $\kappa$ for which it holds
\begin{linenomath*}
\begin{equation}
f(\xi) = (\xi-1)^{\delta_1}(\xi+1)^{\delta_{-1}} \kappa(\xi) + f(1) \frac{1+\psi(\xi)}{2}+ f(-1) \frac{1-\psi(\xi)}{2}.
\end{equation}
\end{linenomath*}
In this subsection we again consider the nonsymmetrical case (\ref{eq2L}-\ref{eq2K}). We denote by $C$ a universal constant, which depends on the data but does not depend on the particular function we consider.

We know that $f$ and hence $\kappa$ are $C^{k-1,\alpha}$ smooth on $\bdr\D\setminus\{-1,1\}$ and on compact subsets of $\bdr\D\setminus\{1,-1\}$ we get a priori estimates on $\kappa$ by expressing it in terms of $f$.
Hence we need a priori estimates on $\kappa$ near points $\pm 1$.
Also, we know from Section \ref{regularity} that if $\kappa$ is continuous on $\overline{\D}$, then both functions belong to $A^{\beta}(\D)$ for 
\begin{linenomath*}
\begin{equation}\label{alpha1-1}
0< \beta < \min\{|\beta_1|, 1-|\beta_1|, |\beta_{-1}|, 1-|\beta_{-1}|\}.
\end{equation}
\end{linenomath*}
Let us fix $0<\beta<\frac{1}{2}-\beta_0$ that we have a priori estimates on function $f$. 

Recall (\ref{rho-wiggle-w}) and that $t^s\kappa(t) = f(t)$.  Hence 
$\widetilde{\rho}_w(\theta,\kappa)$ is a $C^{\beta}$ function with a priori bounds.
As in (\ref{globalkappa}) we can globally write
\begin{linenomath*}
\begin{equation}
{\rm Re}\left(r e^{ i (v+i Hv)} \frac{\partial\kappa}{\partial\theta}\right) = -e^{-u}e^{-(Hv)}\widetilde{\rho}_{\theta}(\theta, \kappa),
\end{equation}
\end{linenomath*}
where $u$ and $v$ are real $C^{\beta}$ functions with a priori bounds.
To get $L^{p'}$ a priori bounds on $\frac{\partial\kappa}{\partial\theta}$ for some $p'>1$ we will get $L^{p'}$ bounds on 
the right-hand side function $\widetilde{\rho}_{\theta} (\theta,\kappa(\theta))$, that it,
on $\widetilde{\rho}_{t} (t,\kappa(t))$ near $t=0$.
Considering (\ref{pt11}-\ref{pt22}-\ref{pt33}) termwise we get
that $t^{\beta_1-1}$, $t^{-\beta_1}$, $\kappa(t) = t^{\beta_1-1} f(t)$, and all terms with function $g$ are $L^{p'}$ bounded for any $p'>1$ such that
\begin{linenomath*}
\begin{equation}
p'(1-\beta_1)<1\ \ \ {\rm and}\ \ \ p'\beta_1<1.
\end{equation}
\end{linenomath*}
Let us consider terms which are bounded by $t^{-\beta_1}|\kappa(t)^2| = t^{\beta_1-2} |f(t)^2|$.
Since we have $\beta\in(0,\frac{1}{2}-\beta_0)$ a priori bounds on $f$, we have
\begin{linenomath*}
\begin{equation}
|f(t)|\le C |t|^{\beta}
\end{equation}
\end{linenomath*}
for some universal constant $C$. Hence
\begin{linenomath*}
\begin{equation}
t^{-\beta_1}|\kappa(t)^2|\le C |t|^{2\beta+\beta_1-2}
\end{equation}
\end{linenomath*}
and this function is in some $L^{p'}$, $p'>1$, if $1<2\beta + \beta_1$. The bound $0<\beta<\frac{1}{2}-\beta_0$ implies that this will be the case for some such $\beta$ if $2\beta_0<\beta_1$. 
Similar argument near $\xi=-1$ gives $2\beta_0<1-|\beta_{-1}|$.

If these two conditions are satisfied, we get $L^{p'}$ a priori estimates on $\widetilde{\rho}_{\theta} (\theta,\kappa(\theta))$ for some $p'>1$.
This implies $C^{\beta'}$ a priori estimate on $\kappa$ for $\beta'<1-\frac{1}{p'}$.

There are natural bounds on $\beta_j$, $j=\pm 1$, in terms of $\beta_0$, that is,
\begin{linenomath*}
\begin{equation}
\frac{1}{2} - \frac{\beta_0}{2}<|\beta_j|<\frac{1}{2}+\frac{\beta_0}{2}.
\end{equation}
\end{linenomath*}
Hence, if $2\beta_0\le \frac{1}{2} - \frac{\beta_0}{2}$ and $\frac{1}{2} + \frac{\beta_0}{2}\le 1-2\beta_0$ both inequalities needed for
$L^{p'}$ a priori estimates will be satisfied. These two inequalities are equivalent to the condition $\beta_0\le\frac{1}{5}$, that is, 
the angle between $w$ and the normal to $\gamma_{\omega_{j}}$
at $w$ is less than $\frac{\pi}{10}$.

\section{Final remarks}

If arc $L$ is the lower semicircle, we can state Theorem \ref{th1} in an equivalent simplified form.
\begin{theorem}\label{th2}
Let $\zp\gamma_{\xi}\zz_{\xi\in\bD\setminus\mathring{L}}$ be a $C^{k}$  $(k\ge 3)$ family of Jordan curves in $\C$ which all contain 
point $0$ in their interiors.
Let Jordan curves $\gamma_{j}$, $j=\pm 1$, be strongly starshaped with respect to $0$ and such that for each 
$w\in \gamma_{\omega_{j}}$ the angle between $w$ and the outer normal to $\gamma_{\omega_{j}}$
at $w$ is less than $\frac{\pi}{10}$.
Let $w_j$, $j=\pm 1$, be the positive intersection of $\gamma_j$ and the real axis with the oriented angle of intersection $\pi\beta_j$, where $\beta_1\in(0,1)$ and $\beta_{-1}\in(-1,0)$. 
Let
\begin{linenomath*}
\begin{equation}
0< \beta < \min\{\beta_1, 1-\beta_1, |\beta_{-1}|, 1-|\beta_{-1}|\}.
\end{equation}
\end{linenomath*}
Then there exists a unique $f\in A^{\beta}(\D)$  with no zeros on $\D$ which solves
(\ref{eq1L}-\ref{eq1K}) for which $f(1)=w_1$ and $f(-1)=w_{-1}$.
\end{theorem}

To prove Theorem \ref{th2} one uses continuity method (see also \cite{For}). The starting boundary value problem (\ref{eq1L}-\ref{eq1K}) can be, using an isotopy from Jordan curves $\zp \gamma_{\xi}\zz_{\xi\bdr\D\setminus\mathring{L}}$ to circles with center at $0$ and fixed radius $R>0$, embedded in a one parameter family of boundary value problems 
which all satisfy assumptions of Theorem \ref{th2}. 
Here, for $t=0$ we have the starting boundary value problem and for $t=1$ circles as the boundary data.

Results in Section \ref{linear-cherepanov} (Proposition \ref{prop2}, Proposition \ref{prop3}) imply that a solution of the boundary value problem (\ref{eq1L}-\ref{eq1K}) for curves 
$\zp\gamma^t_{\xi}\zz_{\xi\bdr\D\setminus\mathring{L}}$
can be locally perturbed into a solution for the nearby perturbed boundary data. Hence the set of parameters $t$ for which there is a solution of  (\ref{eq1L}-\ref{eq1K}) is open. 
On the other hand, a priori estimates from Section \ref{apriori} together with compact embeddings (\ref{compact}) imply that the set of parameters 
$t\in\lbrack 0,1 \rbrack$ for which there is a solution of  (\ref{eq1L}-\ref{eq1K}) is closed. Since there is an obvious solution for the case $t=1$, where all Jordan curves are circles with center at $0$ and fixed radius $R>0$, we get that there is a solution of (\ref{eq1L}-\ref{eq1K}) for $t=0$.

\begin{corollary}\label{cor2}
Let $a_1,\dots, a_n\in\D$ be a finite set of points with given multiplicities.
Then under the assumptions of Theorem \ref{th2}
there exists $\beta\in (0,1)$ and $f\in A^{\beta}(\D)$ which has zeros exactly at points $a_1,\dots, a_n\in\D$ with the given multiplicites
and which solves (\ref{eq2L}-\ref{eq2K}).
\end{corollary}

To prove the corollary we search for solutions $f$ 
of the symmetric problem on the unit disc of the form
\begin{linenomath*}
\begin{equation}
f(z) = \frac{z-a}{1-\overline{a}z}\,  \frac{z-\overline{a}}{1-az}\, \widetilde{f}(z),
\end{equation}
\end{linenomath*}
where $a$ is a point in the upper half-disc. Then $\widetilde{f}$ has to solve a modified problem, where the boundary curves are given by
\begin{linenomath*}
\begin{equation}\label{modified}
\widetilde{\gamma}_{\xi} =  \frac{1-\overline{a}\xi}{\xi-a}\,  \frac{1-a\xi}{\xi-\overline{a}}\, \gamma_{\xi}.
\end{equation}
\end{linenomath*}
Observe that $\widetilde{\gamma}_{\xi} = \gamma_{\xi}$ for $\xi=\pm 1$.

\begin{remark}
{\rm In a similar way one can create a zero at a point $a\in(-1,1)$, that is, on $\mathring{L}$ in the original problem. Jordan curves for the modified problem are
\begin{linenomath*}
\begin{equation}\label{modified}
\widetilde{\gamma}_{\xi} =  \frac{1-a\xi}{\xi-a}\, \gamma_{\xi}.
\end{equation}
\end{linenomath*}
Then $\widetilde{\gamma}_{1}=\gamma_1$ and $\widetilde{\gamma}_{-1} = -\gamma_{-1}$ but conditions of Theorem \ref{th1} are still satisfied.}
\end{remark}

\section*{Acknowledgements}
The author is grateful to the referee for his/her valuable suggestions and comments. 

\section*{Disclosure}
No potential conflict of interest was reported by the author.

\section*{Funding}
The author acknowledges the financial support from the 
Slovenian Research Agency (grants P1-0291, J1-3005, N1-0237 and N1-0137).

\bibliographystyle{plainnat}
\bibliography{Cerne-bib}{}
\end{document}